\documentclass[12pt,a4paper]{article}
\usepackage{amssymb,amsmath}
\def\eqn#1{(\ref{#1})}
\newcommand{\GL}{\mathop{\rm GL}\nolimits}
\newcommand{\IC}{{\mathbb C}}

\newcommand{\II}{{\mathbb I}}

\newcommand{\IP}{{\mathbb P}}

\newcommand{\IT}{{\mathbb T}}
\newcommand{\IZ}{{\mathbb Z}}
\newcommand{\IR}{{\mathbb R}}

\def\abs#1{{\vert#1\vert}}

\def\inf1{{\cal L}^{(1,\infty)}}

\def\ot{\otimes}

\def\pa{\partial}

\def\bra#1{\left\langle #1\right|}
\def\ket#1{\left| #1\right\rangle}
\def\hs#1#2{\left\langle #1,#2\right\rangle}

\def\ca{{\cal A}}
\def\cb{{\cal B}}

\def\ce{{\cal E}}

\def\cm{{\cal M}}

\def\cp{{\cal P}}
\def\cs{{\cal S}}

\def\raw{\rightarrow}

\newbox\ncintdbox \newbox\ncinttbox
\setbox0=\hbox{$-$} \setbox2=\hbox{$\displaystyle\int$}
\setbox\ncintdbox=\hbox{\rlap{\hbox to
\wd2{\hskip-.125em\box2\relax \hfil}}\box0\kern.1em}
\setbox0=\hbox{$\vcenter{\hrule width 4pt}$}
\setbox2=\hbox{$\textstyle\int$}
\setbox\ncinttbox=\hbox{\rlap{\hbox to
\wd2{\hskip-.175em\box2\relax \hfil}}\box0\kern.1em}
\newcommand{\ncint}{\mathop{\mathchoice{\copy\ncintdbox}%
{\copy\ncinttbox}{\copy\ncinttbox}{\copy\ncinttbox}}\nolimits}
\begin{document}

\title{Non-linear $\sigma$-models in noncommutative geometry:  
fields with values in finite spaces}

\author{Ludwik D\c abrowski$\strut^{1}$, Thomas Krajewski$\strut^{2}$, 
Giovanni Landi$\strut^{3}$ \\[15pt] 
$\strut^{1}$ Scuola Internazionale Superiore Studi Avanzati\\
Via Beirut 2-4, I-34014, Trieste, Italy\\
\texttt{dabrow@sissa.it} \\[5pt]
$\strut^{2}$ Centre de Physique Th\'eorique, Campus de Luminy \\
F-13288 Marseille Cedex 9, France\\
and Universit\'e de Provence, France\\
\texttt{Thomas.Krajewski@cpt.univ-mrs.fr} \\[5pt] 
$\strut^{3}$ Dipartimento di Scienze Matematiche, Universit\`a
di Trieste \\ P.le Europa 1, I-34127, Trieste, Italy
\\ and INFN, Sezione di Napoli, Napoli, Italy. \\
\texttt{landi@univ.trieste.it}
}

\maketitle

\begin{abstract}
We study $\sigma$-models on noncommutative spaces, notably on 
noncommutative tori. We
construct instanton solutions carrying a nontrivial topological 
charge $q$ and satisfying a
Belavin-Polyakov bound. The moduli space of these instantons is 
conjectured to consists of  an
ordinary torus endowed with a complex structure times a projective 
space $\IC \IP^{q-1}$.
\end{abstract}

\vfill
Dedicated to A.P. Balachandran on the occasion of his 65th birthday.

\thispagestyle{empty}
\newpage

\section{Introduction}

In \cite{DKL} we have constructed some
noncommutative analogues of two dimensional non-linear 
$\sigma$-models. Since these models
  exhibit rich and easily accessible geometrical structures,  their
noncommutative counterparts are very useful to probe  the interplay
between noncommutative geometry and field theory. We proposed  three classes of
models: an analogue of the Ising model which admits instanton 
solutions, the analogue of
the principal chiral model together with its infinite number of 
conserved currents and the
noncommutative Wess-Zumino-Witten model together with its modified 
conformal invariance.

In this short report we extend the first model by constructing 
instanton solutions for any value of
the topological charge (restricted to be $1$ in \cite{DKL}) and for an
arbitrary complex structure parametrized by a complex number
$\tau \in \IC$, $\Im \tau > 0$.

\section{The general construction}

Ordinary non-linear $\sigma$-models
are field theories whose configuration space consists of maps $X$ 
from the {\em source} space,  a
Riemannian manifold $(\Sigma, g)$, which we assume to be compact and
orientable, to a {\em target} space, an other Riemannian manifold 
$(\cm, G)$. The
corresponding action functional is given, in local coordinates, by
\begin{equation}\label{sigmaaction}
S[X] = {1 \over 2\pi} \int_{\Sigma}
\sqrt{g} ~g^{\mu\nu}\,G_{ij}(X)\partial_{\mu}X^{i}\,\partial_{\nu}X^{j} \, ,
\end{equation} where
as usual $g=\det g_{\mu\nu}$ and $g^{\mu\nu}$ is the inverse of $g_{\mu\nu}$.
The stationary points of this functional are harmonic maps from $\Sigma$ to
$\cm$ and describe minimal surfaces embedded in $\cm$.
Different choices of the source and target spaces lead to different 
field theories,
some of them playing a major role in physics.

When $\Sigma$ is two dimensional, the action $S$ is conformally 
invariant, that is
it is left invariant by any
rescaling of the metric $g\;\rightarrow\;ge^{\sigma}$, where $\sigma$ 
is any map
from $\Sigma$ to $\IR$ . As a consequence, the action only depends
on the conformal class of the metric and may be rewritten using a 
complex structure
on $\Sigma$ as
\begin{equation} S[X] = {i \over \pi}\int_{\Sigma}\,G_{ij}(X)\,\partial
X^{i}\wedge\bar{\partial}X^{j} \, , \end{equation}
where $\partial=\partial_{z}dz$ and
$\bar{\partial}=\partial_{\bar{z}}d\bar{z}$ 
and $z$ is a suitable local complex coordinate.

\bigskip
In order to 
A noncommutative generalization is constructed
by a dualization and reformulation it in terms of the $*$-algebras
$\ca$ and $\cb$ of complex valued smooth functions defined 
respectively on $\Sigma$ and $\cm$. 
Then, embeddings $X$ of $\Sigma$ into $\cm$ 
correspond to $*$-algebra morphisms $\pi_X$ from $\cb$ to $\ca$, 
the correspondence being simply given by the pullback, $f \mapsto \pi_{X}(f)=f\circ X$. 
\\ Now, all this makes perfectly sense also for noncommutative algebras
$\ca$ and $\cb$, and we take as configuration space the space of all 
$*$-algebra morphisms from $\cb$ to $\ca$. 
Both algebras are over $\IC$ and for simplicity 
we take them to be unital. 
The definition of the action functional involves 
noncommutative generalizations of the
conformal and Riemannian geometries. 
According to Connes \cite{book,co1}, the former can
be understood within the framework of positive Hochschild cohomology. 
Indeed, in the commutative two dimensional situation, the tri-linear map 
$\phi: \ca^{\ot 3} \raw \IR$ defined by
\begin{equation}\label{trilinear}
\phi(f_{0},f_{1},f_{2})= {i \over \pi} \int_{\Sigma}f_{0}\partial
f_{1}\wedge\bar{\partial}f_{2}
\end{equation} is an extremal element of the space of positive 
Hochschild cocycles that
belongs to the Hochschild cohomology class of the cyclic cocycle 
$\psi$ defined by
\begin{equation}\label{trilinearbis}
\psi(f_{0},f_{1},f_{2})= {i \over 2 \pi}\int_{\Sigma}f_{0} df_{1} 
\wedge df_{2}.
\end{equation}
Now, expressions
\eqn{trilinear} and \eqn{trilinearbis} still make perfectly sense for 
a general noncommutative algebra $\ca$. One can say that
$\psi$ allows to integrate 2-forms in dimension 2,
so that it is a metric independent
object, whereas $\phi$ defines a suitable scalar product
\begin{equation}\label{3cyco}
\langle a_{0}da_{1},b_{0}db_{1}
\rangle=\phi(b_{0}^{*}a_{0},a_{1},b_{1}^{*})
\end{equation} on the space of 1-forms and thus
depends on the conformal class of the metric. Moreover, this 
scalar product is
positive and invariant with respect to the action of the unitary 
elements of $\ca$
on 1-forms, and its relation to the cyclic cocycle $\psi$ allows to 
prove various 
inequalities involving topological quantities \cite{book}.

We can compose such a cocycle $\phi$ with a morphism
$\pi:\;\cb\rightarrow\ca$ to obtain a positive cocycle on $\cb$ defined by
$\phi_{\pi}=\phi\circ (\pi\ot\pi\ot\pi)$. In order to build an action
functional, which assigns a number to any morphism $\pi$, we have to 
evaluate the
cocycle $\phi_{\pi}$ on a suitably chosen element of $\cb^{\ot 3}$. 
Such an element
provides the noncommutative analogue of the metric on the target, and we
take it as a positive element $G=\sum_{i}b_{0}^{i}\delta 
b_{1}^{i}\delta b_{2}^{i}$
of the space of universal 2-forms $\Omega^{2}(\cb)$. Thus, the quantity
\begin{equation}\label{ncaf}
S[\pi]=\phi_{\pi}(G)
\end{equation}
is well defined and positive. We shall consider it to be a 
noncommutative analogue of
the action functional of the non linear $\sigma$-model.

Clearly, we consider $\pi$
as the dynamical variable (the embedding) whereas $\phi$ (the 
conformal structure
on the source) and $G$ (the metric on the target) are background 
structures that
have been fixed.
Alternatively, one could take only the metric $G$ on the target as a
background field and use the morphism $\pi:\;\cb\rightarrow\ca$ to define the
induced metric $\pi_{*}G$ on the source as
$\pi_{*}G=\sum_{i}\pi(b_{0}^{i})\delta\pi( b_{1}^{i})\delta\pi(
b_{2}^{i})$,
which is obviously a positive universal 2-form on $\ca$. To such an object one
can associate, by means of a variational problem \cite{book,co1}, a 
positive Hochschild
cocycle that stands for the conformal class of the induced metric.

The critical points of the $\sigma$-model corresponding to the action 
functional (\ref{ncaf}) are
noncommutative generalizations of harmonic maps and
describe ``minimally embedded surfaces'' in the noncommutative space 
associated with $\cb$.

\section{Two points as a target space}

The simplest example of a target space is that of a finite space
made of two points $\cm = \{1,2\}$, like in the Ising model. Now, any 
continuous
map from a connected surface to a discrete space is constant and the resulting
(commutative) theory would be trivial. However, this is not the case
if the source is a noncommutative space and one has, 
in general, lots of algebra morphisms.
The algebra of functions over $\cm = \{1,2\}$ is just $\cb = \IC^{2}$ 
and any element $f\in\cm$
is a couple of complex numbers $(f_1,f_2)$ with $f_i=f(i)$, the value 
of $f$ at the point $i$. As a
vector space $\cb$ is generated by the function $e$ defined by 
$e(1)=1, e(2)=0$. Clearly, $e$
is a hermitian projection, $e^2 = e^* = e$, and $\cb$ can be thought 
of as the unital $^*$-algebra
generated by such an element $e$. As a consequence,  any
$*$-algebra morphism $\pi$ from $\cb$ to $\ca$ is given by a 
hermitian projection
$p=\pi(e)$ in $\ca$. The configuration space is then the collection 
of all such projections and our
noncommutative $\sigma$-model provides a dynamics of projections. 
Choosing the metric $G=\delta
e\delta e$ on the space
$\cm$ of two points, the action functional (\ref{ncaf}) simply becomes
\begin{equation} S[p]=\phi(1,p,p),
\end{equation} where $\phi$
is a given Hochschild cocycle corresponding to the conformal structure. \\
As we have already mentioned, from general consideration of 
positivity in Hochschild cohomology this
action  is bounded by a topological term \cite{book}. In the 
following, we shall explicitly prove
this fact when taking the  noncommutative torus as source space.

\section{The noncommutative torus as a source space}\label{se:nct}
We recall the very basic aspects of the
noncommutative torus that we shall need in the following.
The algebra $\ca_\theta$ of smooth functions on the noncommutative 
torus is the unital
$*$-algebra made of power series,
\begin{equation}
a =\sum_{m,n \in \IZ^2} a_{mn} ~(U_1)^m (U_2)^n~,
\end{equation}
with $a_{mn}$ a complex-valued Schwarz function
on $\IZ^2$ that is, the sequence of complex numbers $\{a_{mn} \in 
\IC~, ~ (m,n) \in
\IZ^2 \}$ decreases rapidly at `infinity'. The two unitary elements 
$U_1, U_2$ satisfy the commutation relations
\begin{equation}\label{nct}
U_2 ~U_1 = e^{2\pi i \theta} U_1 ~U_2~.
\end{equation}
There exist on $\ca_\theta$ a 
unique normalized positive definite trace, denoted by the 
integral symbol $\ncint : \ca_\theta \raw \IC$, which is given by
\begin{equation}
\ncint ( \sum_{(m,n) \in \IZ^2} a_{mn} ~(U_1)^m (U_2)^n) := a_{00}~.
\end{equation}
This trace is invariant under the action of the commutative torus 
$\IT^2$ on $\ca_\theta$
whose infinitesimal form is given by two commuting derivations 
$\pa_1, \pa_2$ acting
by\begin{equation}\label{t2act}
\pa_\mu (U_\nu) = 2\pi i ~\delta_\mu^\nu U_\nu~, ~~\mu,\nu = 1, 2~
\end{equation}
and the invariance means just that $\ncint \pa_\mu(a) = 0~, 
~~ \mu = 1, 2~$ for
any $a\in\ca_\theta$.\\
The cyclic 2-cocycle allowing the integration of 2-forms is simply given by
\begin{equation}\label{fc}
\psi(a_{0},a_{1},a_{2})=-{ 1\over 2 \pi i } \ncint \epsilon_{\mu\nu}
a_{0}\partial_{\mu}a_{1}\partial_{\nu}a_{2},
\end{equation}
where $\epsilon_{\mu\nu}$ is the
standard antisymmetric tensor. Its normalization ensures that for any hermitian
projection $p\in\ca_{\theta}$, the quantity $\psi(p,p,p)$ is an 
integer: it is indeed
the index of a Fredholm operator \cite{co}.

The conformal class of a general constant metric is parametrized by a 
complex number
$\tau \in \IC$, $\Im \tau > 0$. Then, up to a conformal factor, the metric is
given by
\begin{equation}\label{met}
g = (g_{\mu\nu}) =
\left(
\begin{array}{cc}
1 & \Re\tau \\
\Re\tau & \abs{\tau}^2
\end{array}
\right)~,
\end{equation}
with inverse given by
$g^{-1} = (g^{\mu\nu}) =
{1 \over (\Im \tau)^2} {\scriptstyle        
  \addtolength{\arraycolsep}{-.5\arraycolsep}
  \renewcommand{\arraystretch}{0.5}
  \left( \begin{array}{cc}
  \scriptstyle \abs{\tau}^2  & \scriptstyle -\Re\tau \\
  \scriptstyle -\Re\tau  & \scriptstyle 1  \end{array} \scriptstyle\right)}$
and $\sqrt{det g} = \Im \tau$.\\
By using the two derivations $\pa_1, \pa_2$ defined in \eqn{t2act} we 
still think
of `the complex torus' $\IT^2$ as acting on the noncommutative torus 
$\ca_\theta$
and construct two associated derivations of
$\ca_\theta$ given by
$\pa_{(\tau)} = {1 \over (\tau - \bar{\tau})} ~(- \bar{\tau} \pa_1 + 
\pa_2)$ and
$\bar{\pa}_{(\tau)} = {1 \over (\tau - \bar{\tau})} ~(\tau \pa_1 - 
\pa_2)$. One easily finds that
$\pa_{(\tau)} \bar{\pa}_{(\tau)} = \bar{\pa}_{(\tau)} \pa_{(\tau)} = 
{1 \over 4} g^{\mu\nu} \pa_\mu
\pa_\nu = {1 \over 4} \Delta$,
and the operator $\Delta = g^{\mu\nu} \pa_\mu \pa_\nu$ is just the 
Laplacian of the
metric \eqn{met}.

By working with the metric (\ref{met}), the positive Hochschild cocycle
$\phi$ associated with the cyclic one (\ref{fc}) will be given by
\begin{equation}\label{phc}
\phi(a_{0},a_{1},a_{2})=
\frac{2}{\pi}\ncint a_{0}\pa_{(\tau)} a_{1}\bar{\pa}_{(\tau)}a_{2}~.
\end{equation}
A construction of the cocycle \eqn{phc} as the conformal class of a general
constant metric on the noncommutative torus can be found in \cite{book,co1}.

\subsection{The action and the field equations}

With $\cp_\theta = Proj(\ca_\theta)$ denoting the collection of all 
projections in the
algebra
$\ca_\theta$, we construct an action functional $S : \cp_\theta \raw \IR^+$ by
\begin{equation}\label{actfun1}
S_{(\tau)}(p) = \phi(1,p,p) = \frac{2}{\pi} \ncint \sqrt{detg} ~ \pa_{(\tau)}
p \bar{\pa}_{(\tau)} p ~.
\end{equation}
The action functional can also be written as
\begin{equation}\label{actfun}
S_{(\tau)}(p) = {1 \over 2 \pi} \ncint \sqrt{detg} ~g^{\mu\nu} \pa_\mu 
p \pa_\nu p =
{1 \over \pi} \ncint \sqrt{detg} ~g^{\mu\nu} p \pa_\mu p \pa_\nu p ~.
\end{equation}
Here the two derivations $\pa_\mu$ are the ones defined in \eqn{t2act}
while the metric $g$ is the one in \eqn{met} which carries also the
dependence on the complex parameter $\tau$.

By taking into account the nonlinear nature of the space 
$\cp_\theta$, one finds that the most
general infinitesimal variation of its elements (the tangent vectors) 
is of the form
$\delta p = (1-p) z p + p z^* (1-p)$, with $z$ an arbitrary elements 
in $\ca_\theta$.
Then, simple algebraic manipulations give the equations of motion,
\begin{equation}\label{eom}
p ~\Delta(p) ~-~ \Delta(p) ~p = 0~,
\end{equation}
where $\Delta = g^{\mu\nu} \pa_\mu \pa_\nu$ is the Laplacian.
These are non linear second order `differential' equations and it is 
not simple to give their solutions in a closed form. 
We shall show that the absolute minima of (\ref{actfun}) 
in a given connected component of $\cp_{\theta}$
actually fulfill first order equations which are easier to solve.

\subsection{Topological charges and self-duality equations}

The cyclic 2-cocycle \eqn{fc} assigns to any projection $p\in\cp_\theta$
a `topological charge' (the first Chern number)
\begin{equation}\label{topcha}
\psi(p) := -{1 \over 2 \pi i} \ncint p \Big[ \pa_1(p) \pa_2(p) - 
\pa_2(p)\pa_1(p) \Big] ~\in
~\IZ~.
\end{equation}
Then, the following inequality holds
\begin{equation}\label{bpbou}
S_{(\tau)}(p) \geq 2 \abs{\psi(p)}~.
\end{equation}
Indeed, due to positivity of the trace $\ncint$ and its cyclic 
properties, we have
that
\begin{equation}\label{ineq}
0 \leq  \ncint \sqrt{det g} ~g^{\mu\nu}
\Big[ \pa_\mu(p) ~p \pm i \epsilon_{\mu}^{~\alpha} \pa_\alpha(p) ~p \Big]^*
\Big[ \pa_\nu(p) ~p \pm i \epsilon_{\nu}^{~\beta} \pa_\beta(p) ~p \Big]~.
\end{equation}
By expanding the LHS and comparing with \eqn{actfun} and \eqn{topcha} we
get inequality \eqn{bpbou}. In \eqn{ineq} the symbol 
$\epsilon_{\mu,\nu}$ stands for the volume form
of the metric $g_{\mu,\nu}$.
The inequality \eqn{bpbou}, which gives a lower bound for the action, is the
analogue of the one for ordinary $\sigma$-models \cite{BePo}.
A similar bound for a model on the fuzzy
$2$-sphere was obtained in \cite{bala}.

It is clear that the equality in \eqn{bpbou} occurs when the projection
$p$ satisfies the following {\it self-dual} or {\it anti-self 
duality} equations
\begin{equation}\label{sd0}
\Big[ \pa_\mu p \pm i \epsilon_{\mu}^{~\alpha} \pa_\alpha p \Big] ~p 
= 0 ~~~~~{\rm
and/or}~~~~~ p \Big[ \pa_\mu p \mp i \epsilon_{\mu}^{~\alpha} 
\pa_\alpha p \Big] = 0
\end{equation}
The self
duality equations \eqn{sd0} can be written as
\begin{equation}\label{sd}
\bar{\pa}_{(\tau)} (p) ~p = 0 ~~~~~{\rm and/or}~~~~~ p ~\pa_{(\tau)} (p) = 0~,
\end{equation}
while the anti-self duality one is
\begin{equation}\label{asd}
\pa_{(\tau)} (p) ~p = 0 ~~~~~{\rm and/or}~~~~~ p ~\bar{\pa}_{(\tau)} (p) = 0~.
\end{equation}
It is straightforward to check that either of the equations 
\eqn{sd} and \eqn{asd}
implies the field equations \eqn{eom} as it should.

\section{The instantons }

The connected components of $\cp_{\theta}$ are parametrized by two 
integers $r$ and
$q$ such that $r+q\theta>0$. \cite{canadian} When $\theta$ is 
irrational, the corresponding
projections have trace $r+q\theta$ and the topological
charge $\psi(p)$ appearing in (\ref{topcha}) is just $q$.

Thus our task is to find projections
that belongs to the previous homotopy classes and satisfy (say) the 
self-duality equation
$\pa(p) ~p = 0$ or, equivalently, $p\pa p=0$.
These equations are still non linear and to solve them 
our next step will be a reduction to a linear problem. 
The key point is to identify the algebra $\ca_\theta$ as the
endomorphism algebra of a suitable bundle and to think of any 
projection in it as an operator on such a bundle. 
The bundle in question will be a projective module of finite
type on a different copy $\ca_\alpha$ of the noncommutative torus, 
the two algebras
$\ca_\theta$ and $\ca_\alpha$ being related by Morita equivalence.

\subsection{The modules}
Let us then consider another copy $\ca_\alpha$ of the noncommutative torus with
generators $Z_1, Z_2$ obeying the relation
$Z_2 Z_1 = e^{2\pi i \alpha} Z_1 Z_2$.
When
$\alpha$ is not rational, every finitely generated projective module 
over the algebra
$\ca_\alpha$ which is not free is isomorphic to a Heisenberg 
module \cite{CR,canadian}. Any such a
module
$\ce_{r,q}$ is characterized by two integers $r,q$ which can be taken 
to be relatively coprime with
$q>0$, or $r=0$ and $q=1$. We shall briefly describe them.
As a vector space
\begin{equation}
\ce_{r,q} = \cs(\IR \times \IZ_q) \simeq \cs(\IR) \otimes \IC^q ~,
\end{equation}
the space of Schwarz functions of one continuous variable $s\in\IR$ 
and a discrete one
$k\in\IZ_q$ (we shall implicitly understand that such a variable is 
defined modulo $q$).
By denoting $\varepsilon = r/q - \alpha$, the space
$\ce_{r,q}$ is made into a {\it right} module over $\ca_\alpha$ by
\begin{equation}\label{rgtmod}
(\xi ~Z_1)(s,k) := \xi(s-\varepsilon, k-r)~, ~~~~~
(\xi ~Z_2)(s,k) := e^{2\pi i (s - k/q)} \xi(s,k)~,
\end{equation}
with $\xi \in \ce_{r,q}$; the relations $Z_2 ~Z_1 = e^{2\pi i \alpha} 
Z_1 ~Z_2$ for the
torus $\ca_\alpha$ are easily verified. \\ On $\ce_{r,q}$ one defines 
an $\ca_\alpha$-valued
hermitian structure
$\hs{\cdot}{\cdot}_{\alpha} : \ce_{r,q} \times \ce_{r,q} \raw 
\ca_{\alpha}$, which is
antilinear in the first factor.

It is proven in \cite{CR} that the endomorphism algebra 
$End_{\ca_\alpha}(\ce_{r,q})$,
which acts on the left on $\ce_{r,q}$, can be identified with another 
copy of the
noncommutative torus $\ca_{\theta}$ where the parameter $\theta$ is 
`uniquely' determined
by $\alpha$ in the following way. Since $r$ and $q$ are coprime, 
there exist integer
numbers $a,b\in\IZ$ such that $ar+bq=1$. Then, the transformed 
parameter is given by
$\theta = (a \alpha + b) / (-q\alpha + r)$.
Notice that, given any two other integers $a',b'\in\IZ$ such that 
$a'r+b'q=1$, one would
find that $\theta' - \theta \in \IZ$ so that $\ca_{\theta'} \simeq 
\ca_{\theta}$.
Thus, we are saying that the algebra
$End_{\ca_\alpha}(\ce_{r,q})$ is generated by two operators $U_1, 
U_2$ acting on the
{\it left} on $\ce_{r,q}$ by
\begin{equation}\label{lftmod}
(U_1 \xi)(s,k) := \xi(s-{1/q}, k-1)~, ~~~~~
(U_2 \xi)(s,k) := e^{2\pi i ~(s / \varepsilon - a k) / q} ~\xi(s,k)~,
\end{equation}
and one verifies that
$U_2 ~U_1 = e^{ 2\pi i \theta} U_1 ~U_2$,
the defining relations of the algebra $\ca_{\theta}$.

The crucial fact is that the $\ca_{\theta}$-$\ca_{\alpha}$-bimodule 
$\ce_{r,q}$ is a {\it Morita
equivalence} between the two algebras $\ca_{\theta}$ and 
$\ca_{\alpha}$: there exists also a
$\ca_{\theta}$-valued hermitian structure on
$\ce_{r,q}$, $\hs{\cdot}{\cdot}_{\theta} : \ce_{r,q} \times \ce_{r,q} \raw
\ca_{\theta}$,
which is compatible with the $\ca_{\alpha}$-valued one
$\hs{\cdot}{\cdot}_{\alpha}$,
\begin{equation}\label{morcom}
\hs{\xi}{\eta}_{\theta} ~\zeta = \xi ~\hs{\eta}{\zeta}_{\alpha}~,
\end{equation}
for all ~$\xi, \eta, \zeta \in \ce_{r,q}$.
The hermitian structure $\hs{\cdot}{\cdot}_{\theta}$ is antilinear in 
the second
factor.

\subsection{The constant curvature connection}
A gauge connection on the right $\ca_{\alpha}$-module $\ce_{r,q}$ is 
given by two
covariant derivative operators
$\nabla_\mu : \ce_{r,q} \raw \ce_{r,q}$, $\mu =1,2$
which satisfy a right Leibniz rule
\begin{equation}\label{leirgt}
\nabla_\mu (\xi ~a) = (\nabla_\mu \xi)a + \xi (\pa_\mu a)~, ~~\mu=1, 2~.
\end{equation}
One also requires compatibility with the $\ca_{\alpha}$-valued 
hermitian structure
\begin{equation}\label{comrgt}
\pa_\mu(\hs{\xi}{\eta}_{\alpha}) = \hs{\nabla_\mu \xi}{\eta}_{\alpha} +
\hs{\xi}{\nabla_\mu \eta}_{\alpha}~, ~~\mu=1, 2~.
\end{equation}
A particular connection on the right $\ca_{\alpha}$-module 
$\ce_{r,q}$ is given by the
operators
\begin{equation}\label{con}
(\nabla_1 \xi)(s,k) := {2\pi i \over \varepsilon} s ~\xi(s,k)~,
~~~(\nabla_2 \xi)(s,k) := {d \xi \over d s}(s,k)~
\end{equation}
(the discrete index $k$ is not touched); this connection is of constant
curvature,
\begin{equation}
F_{1,2} := [\nabla_1, \nabla_2] - \nabla_{[\pa_1, \pa_2]} = -{2 \pi i 
\over \varepsilon}
~\II_{\ce_{r,q}}~,
\end{equation}
with $\II_{\ce_{r,q}}$ the identity operator on $\ce_{r,q}$.

Given any connection $\nabla_\mu$ on $\ce_{r,q}$, one can define 
derivations on the endomorphism
algebra $End_{\ca_\alpha}(\ce_{r,q})$ by commutators: 
$\hat{\delta}_\mu (T) := \nabla_\mu \circ T - T
\circ \nabla_\mu~,~\mu = 1,2$,
for any $T \in End_{\ca_\alpha}(\ce_{r,q})$. Then,
by remembering that $End_{\ca_\alpha}(\ce_{r,q}) \simeq \ca_\theta$, 
one finds that
the derivations $\hat{\delta}_\mu$ on $End_{\ca_\alpha}(\ce_{r,q})$ 
determined by the
particular connection (\ref{con}) are  proportional to the generators 
of the infinitesimal action of
the torus $\IT^2$ on $\ca_\theta$, that is the canonical derivations 
\eqn{t2act}:
$\hat{\delta}_\mu (U_\nu) = {2\pi i \over q \varepsilon}~\delta_\mu^\nu U_\nu~,
~~\mu,\nu = 1, 2$.

The holomorphic and anti-holomorphic connections $\nabla_{(\tau)},
\bar{\nabla}_{(\tau)}$ will be the lift of the derivations 
$\pa_{(\tau)}, \bar{\pa}_{(\tau)}$
  with respect to the connection
\eqn{con}.

\subsection{Instantons from Gaussians}
We shall look for solutions of the self-dual equations \eqn{sd}
of the form
\begin{equation}\label{pro0}
p_\psi := \ket{\psi} \hs{\psi}{\psi}^{-1} \bra{\psi}~,
\end{equation}
with $\ket{\psi}$ a `section of a suitable vector bundle' over the 
noncommutative torus
$\ca_\theta$ and $\hs{\psi}{\psi}$ an invertible element in another 
noncommutative torus
which cannot be $\ca_\theta$ itself but rather is Morita equivalent 
to it, being
indeed $\ca_{\alpha}$. We shall take
$\ket{\psi}$ to be an element of the Schwarz space
$\ce_{r,q}$ on which
$\ca_\theta$ acts on the {\it left} as the endomorphism algebra.
Thus,
$\ce_{r,q}$ will be though of as a {\it right} module over the 
algebra $\ca_{\alpha}$ and
$\hs{\psi}{\psi} = \hs{\psi}{\psi}_{\alpha}$ will be
required to be an invertible element of $\ca_{\alpha}$.

Than, let us suppose that $\ket{\psi} \in \ce_{r,q}$ be such that 
$\hs{\psi}{\psi}_{\alpha}$ is
invertible. A simple computation shows that the projection
$p_\psi := \ket{\psi} (\hs{\psi}{\psi}_{\alpha})^{-1} \bra{\psi}$
is a solution of the self-duality equation \eqn{sd} if and only if the element
$\ket{\psi} \in \ce_{r,q}$ obeys the equation
\begin{equation}\label{sdah}
\bar{\nabla} \psi  - \psi \lambda = 0~,
\end{equation}
with $\lambda$ a suitable element in $\ca_{\alpha}$ and $\bar{\nabla}$ the
anti-holomorphic connection.

\bigskip
For $\lambda \in \IC$ there are simple solutions to \eqn{sdah} given 
by generalized Gaussians
\begin{equation}\label{gau}
\psi_\lambda(s,k) = A_k e^{i \tau \pi s^2 / \varepsilon + 
\lambda(\bar{\tau} - \tau) s} ~.
\end{equation}
Any such a function is an element of $\ce_{r,q}$ since $\Im \tau > 0$.
The vector $A=(A_1, \dots, A_q) \in\IC^{q}$ can be taken in the
complex projective space $\IC \IP^{q-1}$ by removing an
inessential normalization. 
Restrictions on the possible values of the constant parameter $\lambda$
will be discussed presently.

A projection corresponding to a Gaussian 
was constructed in \cite{bo} for
the lowest value of the charge $q=1$ and for $\tau=i$ and the 
invertibility of the corresponding
element $\hs{\psi_\lambda}{\psi_{\lambda}}_{\alpha}$ (now $\alpha=- 
1/\theta$) was
proved for any values of $0<\theta<1$ in \cite{wa}. By using the 
methods of \cite{wa} one should be
able to prove invertibility of 
$\hs{\psi_\lambda}{\psi_{\lambda}}_{\alpha} \in \ca_{\alpha}$ for the
most general situation.

\bigskip

Gauge transformations (in fact complexfied ones) are provided by
invertible elements in $\ca_{\alpha}$ acting on the right on $\ce_{r,q}$,
\begin{equation}\label{gt}
\ce_{r,q} \ni \ket{\psi} \raw \ket{\psi^g} = \ket{\psi} g \in 
\ce_{r,q}~, ~~~\forall
~g\in \GL(\ca_{\alpha})~.
\end{equation}
It is clear that projections of the form \eqn{pro0} are invariant 
under gauge transformations.
Notice that we do not require $g$ to be unitary.

Now, let $\ket{\psi}$ be a solution of \eqn{sdah}, $\bar{\nabla} \psi 
- \psi \lambda = 0$;
and let $g\in \GL(\ca_{\alpha})$. Then a simple computation shows 
that the gauge transformed vector
$\ket{\psi^g}$ will be a solution of an equation of the form 
\eqn{sdah}: $\bar{\nabla} \psi^g  -
\psi^g \lambda_g = 0$ with $\lambda_g$ given by
\begin{equation}\label{pargt}
\lambda_g = g^{-1} \lambda g + g^{-1} \bar{\pa}_{(\tau)} g~.
\end{equation}

The moduli space of Gaussians is simply described. It turns out that 
two Gaussians $\psi_\lambda$ and
$\psi_{\lambda'}$ are related by a gauge transformation if an only if 
there exist  integers $(m,n)\in
\IZ^2$ such that
\begin{equation}
\lambda'=\lambda_g ~, ~~~~~{\rm with}~~~~~ g=(Z_1)^m (Z_2)^n ~.
\end{equation}
Furthermore,
\begin{equation}\label{congt}
\lambda_g - \lambda = {2 \pi i \tau \over {\tau - \bar{\tau}}}
~\left(m - {1 \over \tau} n\right)
\end{equation}
Thus, gauge nonequivalent constant parameters $\lambda$ form a 
complex ordinary torus $\IT_\tau^2$ and
the moduli space of Gaussians \eqn{gau} is $\IC \IP^{q-1} \times \IT_\tau^2$.

\bigskip
It is an open problem to prove whether it is possible to gauge to a 
Gaussian any solution of the
self-duality equation \eqn{sdah}. Equivalently, this problem can be 
stated as follows: given any
$\rho\in\ca_{\alpha}$, there exists an element $g \in 
\GL(\ca_{\alpha})$ such that
$\rho = \lambda  + g^{-1} \bar{\pa}_{(\tau)} g$, with 
$\lambda\in\IC$. Such a gauge transformation can
be found if the deformation parameter $\theta$ is small enough.

\bigskip
As a final remark, me mention that for special values of the 
deformation parameter, $\theta=1/N$, with
$N\in\IZ$, the Gaussian (Boca) projection can be mapped, in the
large torus limit \cite{ks},  to the basic GMS soliton \cite{gms}.

\section*{Acknowledgments}

The work of L.D. and G.L. was supported in part by the 
{\it Progetto INdAM 2003 di ricerca interdisciplinare GNFM-GNSAGA}.

\end{document}